\documentclass[10pt,a4paper,twoside]{article}
\usepackage{amsmath}
\usepackage{amsfonts}
\usepackage{amssymb}
\usepackage{amsthm}
\usepackage{enumerate}
\usepackage{amscd}
\usepackage{amsrefs}
\usepackage{hyperref}
\usepackage{verbatim}

\fontfamily{Liberation Serif}

\hypersetup{citebordercolor={0 1 0}}
\newcommand{\TheoParag}{\newline \hspace*{4 mm}}

\begin{document}
   
    \date{}
    \title{\normalsize \textbf{ON THE GRADED ANNIHILATORS OF RIGHT MODULES OVER THE FROBENIUS SKEW POLYNOMIAL RING}}
    \markboth{\small \textsc{GRADED ANNIHILATORS  }}{\small    \textsc{Ehsan Tavanfar}}         
    \author{\normalsize \textsc{Ehsan Tavanfar}\\ \\ \normalsize \textsc{Dedicated to Professor Hossein Zakeri}\\}
    \maketitle          
    
    \begin{abstract}
       \textsc{Abstract}. Let $R$ be a commutative Noetherian ring of prime characteristic and $M$ be an $x$-divisible right $R[x,f]$-module that is Noetherian as $R$-module. We give an affirmative answer to the question of Sharp and Yoshino in the case where $R$ is semi-local and prove that the set of  graded annihilators of $R[x,f]$-homomorphic images of $M$ is finite. We also give a counterexample  in the	 general case.\newline 
    \end{abstract}     
    \section*{\centering \normalsize \textsc{0. Introduction}}
    \let\thefootnote\relax\footnotetext{2010 Mathematics Subject Classification. 13A35, 13E05, 16S36.}
    \let\thefootnote\relax\footnotetext{Key words and phrases. Commutative Noetherian ring, prime characteristic, skew polynomial ring.}
      \paragraph{} Throughout this note, $R$ is a commutative Noetherian ring of prime characteristic $p>0$ with identity.  Also $f:R\rightarrow R$ is the Frobenius endomorphism which maps an element $r\in R$ to its $p$th power. Let $x$ be an indeterminate. By the Frobenius endomorphism we can endow the free $R$-module generated by $(x^i)_{i\in \mathbb{N}_0}$ with a ring structure such that $xr=f(r)x=r^px$. The result is called the Frobenius skew polynomial ring. We denote it by $R[x,f]$. We can consider $R[x,f]$ as a positively-graded ring $\bigoplus\limits_{i\in\mathbb{N}_{0}}Rx^{n}$.
      
       When $R$ is a local ring with maximal ideal $\mathfrak{m}$, the top local cohomology module of $R$ with respect to $\mathfrak{m}$ has a left $R[x,f]$-module structure. In \cite{1} Rodney Y. Sharp has used this structure to study the existence of the tight closure test element and parameter test element in certain cases.  In order to accomplish these results, he introduced the graded annihilators of left $R[x,f]$-modules. The graded annihilator of a left $R[x,f]$-module $H$ is the largest graded two-sided ideal of $R[x,f]$ contained in $ann_{R[x,f]}H$. The theory of graded annihilators is a novel tool to study the test elements in local rings of prime characteristic: see \cite{1},\cite{2},\cite{3}, for example. One of the key points in this theory is the following finiteness condition (\cite[Corollary 3.11]{1}): if $H$ is an $x$-torsion-free left $R[x,f]$-module which is either Noetherian or Artinian as $R$-module, then the set of  graded annihilators of $R[x,f]$-submodules of $H$ is finite.

       In \cite{4}, Sharp and Yoshino have shown that if $(R,\mathfrak{m})$ is an $F$-finite local ring and $E(R/\mathfrak{m})$ is the injective envelope of the residue field of $R$, then the restriction of the Matlis duality functor $(-)^{\vee}:=Hom_R(-,E(R/\mathfrak{m}))$ to the left $R[x,f]$-modules provides a contravariant functor from the category of left $R[x,f]$-modules to the category of right $R[x,f]$-modules. A similar argument holds for the restriction of the Matlis duality to the right $R[x,f]$-modules.  Also they proved that if $H$ is a left (right) $R[x,f]$-module, then the natural evaluation map $\omega_H:H\rightarrow (H^\vee)^\vee$ is a left (right) $R[x,f]$-homomorphism. By applying these  interesting results, Sharp and Yoshino deduced \cite[Theorem 3.5]{4}: Assume that $R$ is $F$-finite, local and complete. Let $M$ be an $x$-divisible right $R[x,f]$-module which is Noetherian as $R$-module. Then there are only finitely many graded annihilators of $R[x,f]$-homomorphic images of $M$.
       
        This leads to the following question which is proposed in \cite{4}.
    \paragraph{Question.} Assume that $R$ is a commutative Noetherian ring of prime characteristic $p>0$. Let $M$ be an $x$-divisible right $R[x,f]$-module that is Noetherian as an $R$-module. Is the set of  graded annihilators of $R[x,f]$-homomorphic images of $M$ finite?

   Note that this question can be considered as a dual to the finiteness condition mentioned above.  We give an affirmative answer to this question in the case where $R$ is semi-local. Furthermore, we give a counterexample in the general case.
    
    \begin{center}
       \section{\normalsize The Results}
    \end{center} 
    \newtheorem{thm}{Theorem}[section]
    \theoremstyle{Definition}    
    \newtheorem{defi}[thm]{Definition}
    \theoremstyle{Definition and Remark}
    \newtheorem{defi-rem}[thm]{Definitions and Remark}
    \newtheorem{defi-rems}[thm]{Definitions and Remarks}
    \newtheorem{defi-Nots}[thm]{Definition and Notations}
    \theoremstyle{Lemma}
    \newtheorem{lem}[thm]{Lemma}
    \theoremstyle{remark}
    \newtheorem{rem}[thm]{Remark}
    \theoremstyle{Corollary}
    \newtheorem{cor}[thm]{Corollary}
    \newtheorem{exam}[thm]{Counterexample}
    \normalsize  \begin{defi-rems}\hypertarget{1}
 One can easily verify that an ideal $\mathfrak{B}$ of $R[x,f]$ is a two-sided graded ideal if and only if there exists an ascending chain $\{\mathfrak{b_n}\}_{n\in \mathbb{N}_0}$ of ideals of $R$ such that $\mathfrak{B}=\bigoplus\limits_{n\in\mathbb{N}_{0}}\mathfrak{b}_{n}x^{n}$.  
      Assume that $K$ is a right $R[x,f]$-module.

  \begin{enumerate}
    \item[] An $R[x,f]$-submodule $N$ of $K$ is said to be a special submodule if $N$ is of the form $K\mathfrak{B}$ for some two-sided graded ideal $\mathfrak{B}$ of $R[x,f]$. 
  
 \item[] We denote the annihilator of $K$ by $ann(K)_{R[x,f]}$. Note that, 
   \begin{center}
     $ann(K)_{R[x,f]}=\{\alpha\in R[x,f]: k\alpha=0\ \text{for all}\ k\in K\}$,
   \end{center}     
      is a  two-sided ideal of $R[x,f]$. The graded annihilator of $K$ is defined as follows:
     \begin{center}
      $ \text{gr-ann} (K)_{R[x,f]}=\{\sum\limits_{i=0}^{n}r_{i}x^{i}:\ n\in \mathbb{N}_0\ \text{and}\ r_i\in R,\ r_{i}x^{i}\in ann(K)_{R[x,f]}\ \text{for all}\ 0\le i\le n\}$.
   \end{center}
   
     In fact, the graded annihilator of $K$ is the largest graded two-sided ideal of $R[x,f]$ which annihilates $K$.
    
    \item[] We say that $K$ is $x$-divisible if $K=Kx$.
   \end{enumerate}  
  \end{defi-rems}
    \begin{lem}\hypertarget{2}{}
     Let $K$ be a right $R[x,f]$-module. Suppose that $\mathfrak{B}$ and $\mathfrak{B^\prime}$ are  two-sided graded ideals of $R[x,f]$ and that $N$ and $N^\prime$ are  $R[x,f]$-submodules of $K$. 
        
        \begin{enumerate}[i)]
          \item  Whenever $\mathfrak{B}\subseteq \mathfrak{B^\prime}$, $K\mathfrak{B}\subseteq K\mathfrak{B^\prime}$.
          \item  Whenever $N\subseteq N^\prime$, $\text{gr-ann}(K/N)_{R[x,f]}\subseteq \text{gr-ann}(K/N^\prime)_{R[x,f]}$.
          \item  $K(\text{gr-ann}(K/N)_{R[x,f]})\subseteq N$.
          \item  There exists a one-to-one order-preserving correspondence $\Gamma$ between graded annihilators of $R[x,f]$-homomorphic images of $K$ and special submodules of $K$ given by
          \begin{center}
            $\Gamma: K\mathfrak{B}\longmapsto\text{gr-ann}(K/K\mathfrak{B})_{R[x,f]}$. 
          \end{center} 
          The inverse map $\Gamma^{-1}$ is given by
            \begin{center}
              $\Gamma^{-1}:\  \text{gr-ann}(K/N)_{R[x,f]}\longmapsto K(\text{gr-ann}(M/N)_{R[x,f]})$.
            \end{center}
         \end{enumerate}
      \end{lem}
      \begin{proof} Statements (i) and (ii) are obvious. For all  $\alpha \in \text{gr-ann}(K/N)_{R[x,f]}$, we have $K\alpha \subseteq N$ and this   proves $(iii)$. Now, we prove (iv). ‌By parts (ii) and (iii), we have 
        $\text{gr-ann}\Big(K/\big(K(\text{gr-ann}(K/N)_{R[x,f]})\big)\Big)_{R[x,f]}\subseteq\text{gr-ann}(K/N)_{R[x,f]}$.
         For the reverse inclusion suppose that 
           $rx^{n}\in\text{gr-ann}(K/N)_{R[x,f]}$. Then 
             $K(rx^{n})\subseteq K(\text{gr-ann}(K/N)_{R[x,f]})$. 
              So $rx^{n}\in\text{gr-ann}\Big(K/\big(K(\text{gr-ann}(K/N)_{R[x,f]})\big)\Big)_{R[x,f]}$. Therefore 
              \begin{center}
                 $\text{gr-ann}\Big(K/\big(K(\text{gr-ann}(K/N)_{R[x,f]})\big)\Big)_{R[x,f]}=\text{gr-ann}(K/N)_{R[x,f]}$.
              \end{center}   
              On the other hand,  parts (i) and (iii) together with the fact that $\mathfrak{B}\subseteq\text{gr-ann}(K/K\mathfrak{B})_{R[x,f]}$ imply that 
              \begin{center}
                $K\mathfrak{B}=K(\text{gr-ann}(K/K\mathfrak{B})_{R[x,f]})$.
              \end{center}  
              \end{proof}
     \begin{cor} Let $X$ be the the set of all special submodules of a right $R[x,f]$-module $K$. Then $X$ satisfies the ascending chain condition.
     \end{cor}
       \begin{proof} 
       By \cite[Lemma 1.4]{1} (\cite[Lemma 1.4]{1} is deduced by applying the idea of proof of Yoshino's \cite[Corollary (2.7)]{5}), $R[x,f]$ satisfies the ascending chain condition on two-sided graded ideals.  Hence the corollary follows from the previous lemma.
       \end{proof}
     Throughout the rest of this article, we assume that $M$ is an $x$-divisible right $R[x,f]$-module.
     \begin{lem}\hypertarget{3}{}
        Suppose that $0:_RM=\mathfrak{b}$. Then $\mathfrak{b}$ is a radical ideal of $R$ and  $\text{gr-ann}(M)_{R[x,f]}=\mathfrak{b}R[x,f]=\bigoplus\limits_{n\in\mathbb{N}_{0}}\mathfrak{b}x^{n}$.
     \end{lem}
      \begin{proof} There exists an ascending sequence $\{\mathfrak{a}_n\}_{n\in \mathbb{N}_0}$ of ideals of $R$ such that $\mathfrak{a}_0=0:_R M=\mathfrak{b}$ and 
      \begin{center}
        $\text{gr-ann}(M)_{R[x,f]}=\bigoplus\limits_{n\in\mathbb{N}_0}\mathfrak{a}_n x^n$.
      \end{center}        
         Let $r\in\sqrt[]{\mathfrak{a}_0}$. So we have  $r^{p^n} \in \mathfrak{a}_0$ for some $n\in \mathbb{N}_0$, i.e. $Mr^{p^{n}}=0$ for some non-negative integer $n$. But then, 
        $0=(Mr^{p^{n}})x^{n}=M(x^{n}r)=(Mx^{n})r=Mr$. 
         It follows that $\mathfrak{a}_0$ is a radical ideal. Now, it remains to show that $\mathfrak{a}_n\subseteq \mathfrak{a}_0$ for each $n\in \mathbb{N}$. Let $a_n \in \mathfrak{a}_n$, i.e. $M(a_{n}x^{n})=0$. Then $(Ma_{n}^{p^{n}-1})(a_{n}x^{n})=0$ which implies that $(Mx^{n})a_n=0$, and hence $Ma_n=0$. Thus $ a_n\in\mathfrak{a}_{0}$.	\end{proof}
      \begin{defi}
         We say that an ideal $\mathfrak{b}$ of $R$ is a special $M$-ideal if there exists an $R[x,f]$-submodule $N$ of $M$ such that $\bigoplus\limits_{n\in\mathbb{N}_{0}}\mathfrak{b}x^{n}=\text{gr-ann}(M/N)_{R[x,f]}$. 
         
         We will denote the set of all special $M$-ideals by $\mathcal{I}(M)$. Also  we will denote $\mathcal{I}(M)\bigcap Spec(R)$ by $\mathcal{I}^s(M)$. 
      \end{defi}
   
   \begin{cor}
     There exists a one-to-one order-preserving correspondence between $\mathcal{I}(M)$ and the set of special submodules of $M$.
   \end{cor}
   \begin{rem} \hypertarget{10}{} The following hold:
       \begin{enumerate}[i)]
         \item For a family $\{N_i\}_{i\in I}$ of $R[x,f]$-submodules of $M$  we have 
             $\bigcap\limits_{i\in I}(0:_{R}(M/N_{i}))=0:_{R}\big(M/ (\bigcap\limits_{i\in I} N_{i})\big)$.
            So that $\mathcal{I}(M)$  is closed under arbitrary intersection. 	
         \item Special submodules of $M$ are $x$-divisible. To see why this is the case, consider a special submodule $M\mathfrak{B}$ of $M$ where $\mathfrak{B}=\bigoplus\limits_{n\in \mathbb{N}_0}\mathfrak{b}_nx^n$ is a two-sided graded ideal of $R[x,f]$ (Hence $\{\mathfrak{b}_n\}_{n\in\mathbb{N}_0}$ is an ascending chain of ideals of $R$). 
          Let $mb_nx^n\in M\mathfrak{B}$ ($m\in M,\ b_n\in \mathfrak{b}_n,\ n\in \mathbb{N}_0$). Since $M$ is $x$-divisible, we have  $m=m^\prime x$ for some $m^{\prime}\in M$. Hence $mb_{n}x^{n}=(m^{\prime}x)b_{n}x^{n}=m^{\prime}(xb_{n}x^{n})=(m^{\prime}b_{n}^{p}x^{n})x\in(M\mathfrak{B})x$. Therefore, the fact that each element of $M\mathfrak{B}$ can be written  as  a sum of the  elements of the form $mb_nx^n$,  implies that $M\mathfrak{B}\subseteq (M\mathfrak{B})x$.
         \item Assume that $\mathfrak{a},\mathfrak{b}$  are ideals of $R$. Then
           $(M(\mathfrak{a}R[x,f]))(\mathfrak{b}R[x,f])=M((\mathfrak{ab})R[x,f])$.
     \end{enumerate}   
   \end{rem}	 	 
   \begin{thm} \hypertarget{4}{}
     Assume that $N$ is a special submodule of $M$ and  that $\mathfrak{b}$ is the corresponding special $M$-ideal to $N$ (Note that $\mathfrak{b}=0:_R (M/N)$, $N= M(\text{gr-ann}(M/N)_{R[x,f]})=M(\mathfrak{b}R[x,f])$ and that $\mathfrak{b}=\sqrt{\mathfrak{b}}$). Suppose that $\mathfrak{b}=\mathfrak{p}_{1}\bigcap\ldots\bigcap\mathfrak{p}_{n}$ is a minimal primary decomposition of $\mathfrak{b}$ with $n\ge 2$ and that $U$ and $V$ are two non-empty disjoint subsets of $\{1,\ldots,n\}$ such that their union is $\{1,\ldots,n\}$. Let $\mathfrak{a}=\bigcap\limits_{i\in U}\mathfrak{p}_{i}$, $\mathfrak{c}=\bigcap\limits_{i\in V}\mathfrak{p}_{i}$
and $L=M(\mathfrak{a}R[x,f])$. In addition, assume that $\mathfrak{b}^\prime$ is an ideal of $R$ such that $N\subseteq M(\mathfrak{b}^{\prime}R[x,f])$. Then
     \begin{enumerate}[i)]
        \item $\text{gr-ann}(M(\mathfrak{b}^{\prime}R[x,f])/N)_{R[x,f]}=(\mathfrak{b}:_{R}\mathfrak{b^{\prime}})R[x,f]$.
        \item $\text{gr-ann}(L/N)_{R[x,f]}=\mathfrak{c}R[x,f]$. In particular, $\mathfrak{c}\in \mathcal{I}(L)$.
        \item The corresponding $M$-special ideal to $L$ is $\mathfrak{a}$. In particular, $\mathfrak{a}\in \mathcal{I}(M)$.
     \end{enumerate}
   \end{thm}
     \begin{proof}\begin{enumerate}[i)] 
     \item   Since $M(\mathfrak{b}^\prime R[x,f])$ is an $x$-divisible right $R[x,f]$-module, we only need to prove that 
     \begin{center}
       $0:_R M(\mathfrak{b}^\prime R[x,f])/N=\mathfrak{b}:_R\mathfrak{b}^\prime$.
     \end{center}       
        It is easy to verify that 
       $\mathfrak{b}:_R\mathfrak{b}^\prime\subseteq 0:_RM(\mathfrak{b}^\prime R[x,f])/N$.
        For the reverse inclusion, let  $m\in M$,  $r\in 0:_R M(\mathfrak{b}^\prime R[x,f])/N$ and $t\in0:_{R} M/M(\mathfrak{b}^\prime R[x,f])$.    Then $mt\in M(\mathfrak{b}^\prime R[x,f])$ and $(mt)r\in N$. So $tr\in 0:_RM/N$. This implies that $(0:_R\ M/M(\mathfrak{b}^\prime R[x,f]))r\subseteq 0:_RM/N$. But,
        \begin{center}
          $\mathfrak{b}^\prime\subseteq 0:_RM/M(\mathfrak{b}^\prime R[x,f])$.
        \end{center}          
           Hence $\mathfrak{b}^\prime r\subseteq 0:_RM/N=\mathfrak{b}$, i.e. $r\in (\mathfrak{b}:_R\mathfrak{b}^\prime)$.
     \item By \cite[Lemma 3.5(ii)]{1}, $\mathfrak{b}:_R\mathfrak{a}=\mathfrak{c}$.
     \item Suppose that $\text{gr-ann}(M/L)_{R[x,f]}=\mathfrak{d}R[x,f]$. If we prove that $\mathfrak{d}\subseteq \mathfrak{a}$, then we are done. Note that $\mathfrak{d}(0:_{R}L/N)\subseteq0:_{R}M/N=\mathfrak{b}$ and $0:_RL/N=\mathfrak{c}$. Thus  $\mathfrak{d}\mathfrak{c}\subseteq \mathfrak{b}$. Since $\mathfrak{p_{1}}\bigcap\ldots\bigcap\mathfrak{p_{n}}$ is the minimal primary decomposition of $\mathfrak{b}$,  we have $\mathfrak{c}\nsubseteq\mathfrak{p}_{i}$ for each $i\in U$. But $\mathfrak{d}\mathfrak{c}\subseteq\mathfrak{b}\subseteq\mathfrak{p}_{i}$
for each $i\in U$. Therefore $\mathfrak{d}\subseteq\bigcap\limits_{i\in U}\mathfrak{p}_{i}=\mathfrak{a}$. 
  \end{enumerate} \end{proof}
   \begin{cor}\hypertarget{5}{}
    $\mathcal{I}(M)\backslash {\{R\}}$ is precisely the set of all finite intersections of members of $\mathcal{I}^s(M)$.
   \end{cor}
   \begin{thm} \hypertarget{8}{}
     Suppose that $M$ is Noetherian as $R$-module. Then 
     \begin{center} 
       $\mathcal{I}^s(M)=\{\mathfrak{p}\in Spec(R):\ \mathfrak{p}\in  Ass_R(M/M(\mathfrak{p}R[x,f]))\}$.
     \end{center} 
     In particular, $\mathcal{I}^s(M)\subseteq Supp_R(M)$.
   \end{thm}
    \begin{proof} Note that $\mathfrak{p}\in Ass_R(M\slash M(\mathfrak{p}R[x,f]))$ if and only if $\mathfrak{p}=0:_R M\slash M(pR[x,f])$, and the latter is equivalent to that  $\mathfrak{p}\in \mathcal{I}^s(M)$. On the other hand, if $\mathfrak{p}\in Ass_R(M\slash M(\mathfrak{p}R[x,f]))$, then we have $\mathfrak{p}=0:_{R}(m+M(\mathfrak{p}R[x,f]))$ for some $m\in M$ . Hence
     
     \begin{center}
       $\mathfrak{p}\subseteq0:_{R}(m+M\mathfrak{p})\subseteq0:_{R}(m+M(\mathfrak{p}R[x,f]))=\mathfrak{p}$.
     \end{center}  
   Therefore $\mathfrak{p}\in Ass_{R}(M/M\mathfrak{p})$. Thus 
     $\mathcal{I}^s(M)\subseteq \{\mathfrak{p}\in Spec(R):\ \mathfrak{p}\in Ass_R(M/M\mathfrak{p})\}\subseteq Supp_R(M)$.
  \end{proof}

    \begin{thm}\hypertarget{6}{}
      Suppose that  $L:=M(\mathfrak{a} R[x,f])$ for some ideal $\mathfrak{a}\in \mathcal{I}(M)\backslash\{0:_R M\}$. Then
        \begin{center}
          $\mathcal{I}^{s}(L)=\{\mathfrak{p}\in\mathcal{I}^{s}(M):\ \mathfrak{a}\nsubseteq\mathfrak{p}\}$.
        \end{center}
    \end{thm}
     \begin{proof} Let $\mathfrak{p}\in\mathcal{I}^{s}(L)$, i.e. $\mathfrak{p}=0:_{R}L/L(\mathfrak{p}R[x,f])$. Let $b\in 0:_R M/M(\mathfrak{p}R[x,f])$  and $m(ax^{i})\in L$ ($m\in M,\ a\in \mathfrak{a}$ and $ \ i\in \mathbb{N}_0$). By  \hyperlink{10}{1.7(iii)}, we have
      \begin{center}
        $(max^{i})b=(mb^{p^{i}})ax^{i}\in(M(\mathfrak{p}R[x,f]))(\mathfrak{a}R[x,f])=(M(\mathfrak{a}R[x,f]))(\mathfrak{p}R[x,f])=L(\mathfrak{p}R[x,f])$.
      \end{center}
      It follows that 
      $0:_R M/M(\mathfrak{p}R[x,f])\subseteq0:_{R} L/L(\mathfrak{p}R[x,f]) =\mathfrak{p}$. Thus, $0:_R M/M(\mathfrak{p}R[x,f])=\mathfrak{p}$, i.e. $\mathfrak{p}\in \mathcal{I}^s(M)$.
      
      We have 
        $L(\mathfrak{p}R[x,f])=M((\mathfrak{ap})R[x,f])=M(\mathfrak{c}R[x,f])$
      for some radical ideal $\mathfrak{c}$ such that 
      \begin{center}
        $\mathfrak{c}=0:_{R}M/M((\mathfrak{ap})R[x,f])$.
      \end{center}        
         So if $\mathfrak{a}\subseteq \mathfrak{p}$, then 
$\mathfrak{a}^{2}\subseteq\mathfrak{ap}\subseteq\mathfrak{c}\subseteq\mathfrak{a}$, hence $\mathfrak{c}=\mathfrak{a}$ and $L(\mathfrak{p}R[x,f])=L$ which contradicts  the choice of $\mathfrak{p}$. It follows that $\mathfrak{p}\in \{\mathfrak{p}^\prime\in \mathcal{I}^s(M):\   \mathfrak{a}\nsubseteq \mathfrak{p}^\prime\}$.
      
       Conversely,  let $\mathfrak{p}\in \{\mathfrak{p}^\prime\in\mathcal{I}^{s}(M):\ \mathfrak{a}\nsubseteq\mathfrak{p}^\prime\}$. Note that, in particular,
 $\mathfrak{p}:_{R}(\mathfrak{a}+\mathfrak{p})=\mathfrak{p}$. 
 Hence by \hyperlink{4}{Theorem 1.8(i)},
$\mathfrak{p}=0:_{R}M((\mathfrak{a}+\mathfrak{p})R[x,f])/M(\mathfrak{p}R[x,f])$.

 For each 
 $r\in 0:_R L/L(\mathfrak{p}R[x,f]) $ and $m((a+p)x^{i})\in M((\mathfrak{a}+\mathfrak{p})R[x,f])$ ($m\in M,\ a\in\mathfrak{a},\ p\in\mathfrak{p}$ and $i\in \mathbb{N}_0$) we have
 \begin{center} 
      $(m((a+p)x^{i}))r=(m(ax^{i}))r+(m(px^{i}))r\in M(\mathfrak{p}R[x,f])$.
 \end{center}     
     It follows that $0:_{R} L/L(\mathfrak{p}R[x,f]) \subseteq 0:_{R}M((\mathfrak{a}+\mathfrak{p})R[x,f])/M(\mathfrak{p}R[x,f])=\mathfrak{p}$. Thus $p\in \mathcal{I}^s(L)$. \end{proof}
     \begin{thm} \hypertarget{7}{}
       Assume that $R$ is a semi-local ring and that $M$ is Noetherian as $R$-module. Then there are only finitely many graded annihilators of $R[x,f]$-homomorphic images of $M$.
     \end{thm}
      \begin{proof} 
      Firstly, we claim that the set of maximal member of $\mathcal{I}^s(M)$ is finite. Suppose the  contrary,  and  let $\{\mathfrak{p}_i\}_{i\in \mathbb{N}}$ be a countable (infinite) subset of the maximal members of $\mathcal{I}^s(M)$. ‌By our hypothesis, 
            \begin{center}
         $M((\bigcap\limits_{i=1}^{n}\mathfrak{p_{i})}+\mathfrak{p}_{n+1})R[x,f])=M$
      \end{center}   
          for all $n\ge 1$. So using the exact sequence 
     \begin{center}
       
$0\longrightarrow\ M/\bigg(\big(\bigcap\limits_{i=1}^{n}M\mathfrak{(p_{i}}R[x,f])\big)\bigcap M(\mathfrak{p}_{n+1}R[x,f]) \bigg)\longrightarrow \big(M/(\bigcap\limits _{i=1}^{n}M\mathfrak{(p_{i}}R[x,f])) \big)\bigoplus \big(M/M(\mathfrak{p}_{n+1}R[x,f])\big) \longrightarrow  M/\bigg(\big(\bigcap\limits _{i=1}^{n}M\mathfrak{(p_{i}}R[x,f])\big)+M(\mathfrak{p}_{n+1}R[x,f])\bigg) \longrightarrow0$
     \end{center}
     and by induction on $n$, we  deduce that $M/\big(\bigcap\limits _{i=1}^{n}M\mathfrak{(p_{i}}R[x,f])\big)\cong\bigoplus\limits_{i=1}^{n}M/M(\mathfrak{p_{i}}R[x,f])$ for each $n\in \mathbb{N}$.  \TheoParag  We denote the jacobson radical of $R$ by $\mathfrak{a}$. Furthermore, if an $R$-module $L$ is Noetherian, then we use $\nu(L)$ to denote the length of the $R$-module $L/\mathfrak{a}L$. By the previous paragraph, there exists an epimorphism $M\rightarrow \bigoplus\limits_{i=1}^nM/M(\mathfrak{p}_iR[x,f])$ for each $n\in \mathbb{N}$. Therefore $\nu(M)\ge \sum\limits_{i=1}^n\nu\big(M/M(\mathfrak{p}_iR[x,f])\big)$  for all $n\in \mathbb{N}$. Thus there must exist a natural number $i$ such that $\nu\big(M/M(\mathfrak{p}_iR[x,f])\big)=0$, hence $M/M(\mathfrak{p}_iR[x,f])=0$ (by Nakayama's lemma). This is a contradiction. Therefore the set of  maximal members of $\mathcal{I}^s(M)$ is finite.

          By \hyperlink{5}{Corollary 1.9}, we only need to prove that $\mathcal{I}^s(M)$ is finite. Suppose the contrary and let $dim(R)=t$. By the previous paragraph $\mathcal{I}^s(M)$ has only finitely many maximal members. If $t=0$, then each member of $\mathcal{I}^s(M)$ is maximal which is a contradiction. So $t>0$. Let $M_0=M$ and $\mathfrak{a}_1$ be the intersection of the maximal members of $\mathcal{I}^s(M)$. By \hyperlink{10}{Remark 1.7(i)}, $\mathfrak{a}_1\in \mathcal{I}(M)$. Note that $0:_RM\neq \mathfrak{a}_1$, otherwise $\mathcal{I}^s(M)$ is finite (Assume that $\mathfrak{q}\in  \mathcal{I}^s(M)$ and  that $\{\mathfrak{p}_i\}_{i=1}^n$ are maximal members of $\mathcal{I}^s(M)$. If $0:_R M = \mathfrak{a}_1$, then we have
         $\mathfrak{q}\supseteq0:_{R}M=\mathfrak{a}_1=\bigcap\limits_{i=1}^{n}\mathfrak{p}_{i}$.
       Hence $\mathfrak{q}=\mathfrak{p}_i$ for some $i$ by the maximality).
       \TheoParag Let $M_1=M(\mathfrak{a}_1R[x,f])$. By \hyperlink{6}{Theorem 1.11} 
       \begin{center}
         $\mathcal{I}^s(M_1)=\mathcal{I}^s(M)\backslash\{\text{maximal members of }\mathcal{I}^s(M)\}$
       \end{center}
        and $\mathcal{I}^s(M_1)$ is infinite. Assume that, by induction, we have constructed an $x$-divisible $R[x,f]$-submodule $M_t$ of $M$ such that 
        \begin{center}
          $\mathcal{I}^s(M_t)=\mathcal{I}^s(M_{t-1})\backslash\{\text{maximal members of }\mathcal{I}^s(M_{t-1})\}$
        \end{center}         
           and  $\mathcal{I}^s(M_t)$ is infinite. As before, let $\mathfrak{a}_{t+1}$ be the intersection of maximal members of $\mathcal{I}^s(M_t)$ and $M_{t+1}=M_t(\mathfrak{a}_{t+1}R[x,f])$. Note that again $\mathfrak{a}_{t+1}\neq 0:_R M_t$. Therefore $\mathcal{I}^s(M_{t+1})$ is infinite (In particular, it is non-empty). Let $\mathfrak{p}_{t+1}\in \mathcal{I}^s(M_{t+1})$. There exists a maximal member of $\mathcal{I}^s(M_t)$, say $\mathfrak{p}_t$ such that $\mathfrak{p}_{t+1}\subsetneq \mathfrak{p}_t$. Proceeding in this way we can make an ascending chain of prime ideals of $R$ of length $t+1$. This is a contradiction.
          \end{proof}
       \begin{cor}\hypertarget{14} Suppose that $M$ is Noetherian as $R$-module and that
       \begin{center}
          $\mathcal{A}:=\{\mathfrak{m}\in Max(R):\text{ there\  exists\  a\  maximal member }\ \mathfrak{p}\ \in \mathcal{I}^s(M)\   such\  that\  \mathfrak{p}\subseteq \mathfrak{m}\}$
       \end{center}
        is a finite set. Then the set of  graded annihilators of $R[x,f]$-homomorphic images of $M$ is finite.\newline In particular, if the set of maximal members of $\mathcal{I}^s(M)$ is finite, then there are only finitely many graded annihilators of $R[x,f]$-homomorphic images of $M$.
       \end{cor}
       \begin{proof} Suppose the contrary and assume that $\mathcal{I}^{s}(M)$ is not finite. By our hypothesis,  there exists a maximal ideal $\mathfrak{m}$ of  $R$ such that there are infinitely many members of $\mathcal{I}^{s}(M)$
contained in $\mathfrak{m}$. By \cite[Lemma 3.3]{4}, $M_{\mathfrak{m}}$ has an
$x$-divisible right $R_{\mathfrak{m}}[x,f]$-module structure such
that $\frac{m}{s}x=\frac{ms^{p-1}x}{s}$ for all $m\in M$ and $s\in R\backslash\mathfrak{m}$. Since $(M(\mathfrak{q}R[x,f]))_{\mathfrak{m}}=M_{\mathfrak{m}}((\mathfrak{q}R_{\mathfrak{m}})R_{\mathfrak{m}}[x,f])$, by applying \hyperlink{8}{Theorem 1.10}, we  deduce that,
  \begin{center}
    $\mathcal{I}_{R_{\mathfrak{m}}}^{s}(M_{\mathfrak{m}})=\{\mathfrak{q}R_{\mathfrak{m}}:\ \mathfrak{q}\in\mathcal{I}^{s}(M),\       \mathfrak{q}\subseteq\mathfrak{m}\}$.
  \end{center}
  But then $\mathcal{I}_{R_{\mathfrak{m}}}^{s}(M_{\mathfrak{m}})$ is infinite which contradicts the
previous theorem. 
 
  Similarly, the second assertion can be proved by localization at maximal members of $\mathcal{I}^s(M)$. \end{proof}

  The following example shows that \hyperlink{7}{Theorem 1.12} does not hold without the assumption for the ring to be semi-local.
\begin{exam} Consider the polynomial ring
\begin{center}
  $A:=K[\theta_1,\theta_2,\ldots,\theta_n,\ldots]$
\end{center}  
    over a perfect field $K$ of prime characteristic $p>0$, in a countably infinite set of indeterminates. Let $\mathfrak{p}_i:=(\theta_i)$ for each $i\in \mathbb{N}$. Thus $\{\mathfrak{p}_{i}\}_{i\in\mathbb{N}}$ is a family of prime ideals. Therefore we can consider the multiplicative closed subset $S:=A\backslash (\bigcup \limits_{i\in \mathbb{N}}\mathfrak{p}_i)$ of $A$. It is easy to verify that $R:=S^{-1}A$ is a Noetherian ring of dimension 1. We are going to prove that $R$ has an $x$-divisible right $R[x,f]$-module structure and that $S^{-1}\mathfrak{p}_i	\in \mathcal{I}^s(R)$ for each $i\in \mathbb{N}$.  
    \TheoParag The Frobenius map can be used to endow $A$ with a left $A$-module structure, denoted by $F(A)$, such that $b.a=b^p a$ for all $a,b\in A$. We define the map $h:F(A)\rightarrow A$  by 	
    \begin{center}
      $h(\sum\limits_{i=0}^{t}k_{i}\theta_{1}^{\alpha_{i,1}}\ldots \theta_{n_{i}}^{\alpha_{i,n_{i}}})=\sum\limits_{\substack{i=0\\ \alpha_{i,1},\ldots,\alpha_{i,n_{i}}\text{are}\\ \text{multiple of }p}
}^{t}k_{i}^{1/p}\theta_{1}^{\alpha_{i,1}/p}\ldots\theta_{n_{i}}^{\alpha_{i,n_{i}}/p}$,
    \end{center}
    for all 
    $t\in \mathbb{N}_0$, $n_0,\ldots,n_t\in \mathbb{N}_0$, $k_0,\ldots,k_t\in K$ and $\alpha_{0,1},\ldots,\alpha_{0,n_1},\ldots,\alpha_{t,1},\ldots,\alpha_{t,n_t}\in \mathbb{N}_0$.
     It is straightforward to check that $h$ is an $A$-homomorphism which splits the exact sequence 
     \begin{center}
       $0\rightarrow A\overset{f}{\rightarrow} F(A)\rightarrow\frac{F(A)}{A}\rightarrow0$. 
     \end{center}       
       Therefore there exists an $A^p$-epimorphism $\pi:A\rightarrow A^p$ such that $\pi(a)=h(a)^p$ for all $a\in A$. Hence by \cite[Example 1.9(ii)]{4}, $A$ is an $x$-divisible right $A[x,f]$-module such that $ax=\pi(a)^{\frac{1}{p}}=h(a)$ for all $a\in A$. So By \cite[Lemma 3.3]{4}, $R$ has a natural $x$-divisible right $R[x,f]$-module structure.
     \TheoParag It is easy to verify that $\mathfrak{p}_i = A(\mathfrak{p}_i A[x,f])$ for each $i\in \mathbb{N}_0$. So 
       $\mathfrak{p}_i\in Ass_A(A/A(\mathfrak{p}_i A[x,f]))$.
      Thus 
      \begin{center}
        $S^{-1}\mathfrak{p}_i\in Ass_R\big(R/S^{-1}(A(\mathfrak{p}_i A[x,f]))\big)$. 
      \end{center}        
        But 
        $S^{-1}(A(\mathfrak{p}_{i}A[x,f]))=S^{-1}A(S^{-1}\mathfrak{p}_{i}(S^{-1}A[x,f]))=R(S^{-1}\mathfrak{p}_{i}R[x,f])$.
        Hence 
        \begin{center}
          $S^{-1}\mathfrak{p}_{i}\in Ass_{R}(R/R(S^{-1}\mathfrak{p}_{i}R[x,f]))$. 
        \end{center}          
          Therefore $S^{-1}\mathfrak{p}_{i}\in\mathcal{I}^{s}(R)$ by \hyperlink{8}{Theorem 1.10}. This implies that $\mathcal{I}^s(R)$ is infinite. \\ \\
     (The statements of \cite[Example 1.9(ii)]{4} and \cite[Lemma 3.3]{4} still hold if we drop the assumption that $R$ is Noetherian.)
\end{exam}

  \begin{thm} 
     Assume that $M$ is Artinian as $R$-module. Then the set of graded annihilators of $R[x,f]$-homomorphic images of $M$ is finite.
  \end{thm}
  \begin{proof} Suppose, to the contrary, that $\mathcal{I}^s(M)$ is infinite. Let $\mathfrak{p}$ be a maximal member of $\mathcal{I}^s(M)$ and $M_{1}=M(\mathfrak{p}R[x,f])$. Note that $M_{1}\subsetneq M$. By \hyperlink{6}{Theorem 1.11}, $\mathcal{I}^{s}(M_{1})$ is infinite. Let $M_0=M$ and $t\ge1$. Assume,  inductively, that we  have constructed an x-divisible $R[x,f]$-submodule $M_{t}$ of $M$ such that $M_{t}\subsetneq M_{t-1}$ and $\mathcal{I}^{s}(M_{t})$ is infinite.  Let $\mathfrak{p}_{t}$ be a maximal member 
of $\mathcal{I}^{s}(M_{t})$ and  
$M_{t+1}=M_{t}(\mathfrak{p}_{t}R[x,f])$ again. So $M_{t+1}\subsetneq M_{t}$
and $\mathcal{I}^{s}(M_{t+1})$ is infinite. Hence by
induction, we have constructed a strictly descending chain of $R$-submodules
of $M$ which is a contradiction. \end{proof}

\paragraph{Acknowledgements.} I would like to thank and express my deepest gratitude to my MSc supervisor, Professor Hossein Zakeri, for his valuable comments, patience and encouragement during preparation of this paper and my thesis.

\begin{center}\section*{\normalsize \textsc{References}}\end{center}

\small \textsc{Faculty of Mathematical Sciences and Computer, Tarbiat Moallem University, 599 Taleghani Avenue, Tehran 15618, IRAN.} \newline \\
\small \textsc{School of Mathematics, Institute for Research in Fundamental Sciences (IPM),
P.O. Box: 19395-5746, Tehran, IRAN.}\\ \\
E-mail address: \href{mailto:tavanfar@gmail.com}{tavanfar@gmail.com} \newline

\end{document}